\documentclass[10pt]{IEEEtran} %\onecolumn
\usepackage{amsfonts}       % blackboard math symbols
\usepackage{amsmath,amssymb,amsthm, bm, bbm,dsfont}
\usepackage{algorithm}
\usepackage{array}
\usepackage{mdwmath}
\usepackage{url}
\usepackage{multirow}

\usepackage{graphicx}
\usepackage{caption}
\usepackage{subcaption}

\usepackage{tikz}
\usepackage{algpseudocode}
\usepackage{pgfplots}
\newtheorem{theorem}{Theorem}[section]

\newtheorem{fact}[theorem]{Fact}

\newtheorem{remark}[theorem]{Remark}

\newcommand{\bproof}{ \begin{IEEEproof} }
\newcommand{\eproof}{ \end{IEEEproof} }
\newcommand{\beqno}{ \begin{equation*} }
\newcommand{\eeqno}{ \end{equation*} }
\newcommand{\beqa}{\begin{eqnarray*} }
\newcommand{\eeqa}{\end{eqnarray*} }
\newcommand{\beq}{ \begin{equation} }
\newcommand{\eeq}{ \end{equation} }

\renewcommand{\a}{\bm{a}}

\newcommand{\x}{\bm{x}}

\newcommand{\y}{\bm{y}}
\newcommand{\w}{\bm{w}}

\newcommand{\X}{\bm{X}}

\newcommand{\I}{\bm{I}}

\newcommand{\n}{\mathcal{N}}
\newcommand{\E}{\mathbb{E}}

\newcommand{\z}{\bm{z}}

\newcommand{\iidsim}{\stackrel{\mathrm{iid}}{\thicksim }}

\newcommand{\A}{\bm{W}}%{\bm{A}}

\newcommand{\Lambar}{\bar{\bm\Lambda}}
\newcommand{\lammax}{\bar\lambda_{\max}} %{\lambda_{mx}}

\newcommand{\N}{\mathcal{N}}

\newcommand{\W}{\bm{D}} %{\bm{W}}
\renewcommand{\Re}{\mathbb{R}}

\begin{document}
%
%To Do for JP:
%- discuss and cite Zheng and Lafferty. they need total measurements > nr^3 \log n kappa^2. Similar to what we need in the deterministic case. But they need GOE data.
%- add AltMinPhase. compare with idea of Rev 2. compare with using the same a_i's for all x_k's.
%n=1000 sim.
%- estimate r

%\title{Low Rank Matrix Recovery from Column-wise Phaseless Measurements} , Yonina C. Eldar Y. C. Eldar is with Technion, Haifa, Israel.
\title{A Simple Generalization of a Result for Random Matrices with Independent Sub-Gaussian Rows} %Matrices with Independent Sub-Gaussian Rows} Simple  Vershynin's
\author{Namrata Vaswani and Seyedehsara Nayer
\thanks{N. Vaswani and S. Nayer are with the Iowa State University, Ames, IA, USA. Email: namrata@iastate.edu.
}
}

\maketitle

%Still need: explicit algorithm; add video pictures.

%e initialization step of AltMinTrunc is able to get a solution within a given error radius of the true one uses ideas from two recent approaches to phase retrieval for a single vector to solve the current problem. borrows some ideas from two recent phase retrieval works. It
%high probability lower bounds on the number of samples needed by the initialization step of AltMinTrunc (which is also the hardest to analyze typically)
%get within an $\varepsilon$ ball of the true
\begin{abstract} %magnitude-only ((``low rank phase retrieval").
In this short note, we give a very simple but useful generalization of a result of Vershynin (Theorem 5.39 of \cite{vershynin}) for a random matrix with independent sub-Gaussian rows.
We also explain with an example where our generalization is useful.%The proof follow using the strategy outlined for proving Vershynin's result.  $\A$,
%Our result does two extra things. First, it bounds the deviation of $\A$ from its expected value, $\E[A]$, even when the different rows of $\A$ do not have the same second moment matrix. Second, it states a separate result that bounds $\|(\A - \E[A]) \x\|_2^2$ for a given vector $\x$. This bound clearly holds with much higher probability than the bound on $\|\A - \E[A]\|_2$.
%while the result of Vershynin is not.
%Theorem 5.39 of Vershynin's paper on non-asymptotic results for random matrices
\end{abstract}

\section{Introduction}
\label{intro}

In this note, we obtain a generalization of a result of Vershynin, Theorem 5.39 of \cite{vershynin}. This result bounds the minimum and maximum singular values of an $N \times n$ matrix $\A$ with mutually independent, sub-Gaussian, and isotropic rows.
%
%For other $l_p$ norms, we use $\|.\|_p$.
We use $\|.\|$ to denote the $l_2$ norm of a vector or the induced $l_2$ norm of a matrix,  and we use $'$ to denote matrix or vector transpose. Let $\A = [\w_1, \w_2, \dots \w_N]'$. Thus, $\w_j$ is its $j$-th row. As explained in \cite[Section 5.2]{vershynin}, ``isotropic" means that $\E[\w_j \w_j{}' ] = \I$ where $\I$ is the identity matrix. In Remark 5.40 of \cite{vershynin}, this result is generalized to the case where the rows $\w_j$ are not isotropic but have the same second moment matrix, $\E[\w_j \w_j{}' ]$ for all the $N$ rows.
As explained in \cite{vershynin}, a sub-Gaussian random variable (r.v.), $x$, is one for which the following holds: there exists a constant $K_g$ such that
$
\E[|x|^p]^{1/p} \le K_g \sqrt{p}
$
for all integers $p \ge 1$. The smallest such $K_g$ is referred to as the sub-Gaussian norm of $x$, denoted $\|x\|_{\varphi_2}$. Thus,
$
\|x\|_{\varphi_2} = \sup_{p \ge 1} p^{-1/2} \E[|x|^p]^{1/p}.
$
A sub-Gaussian random vector, $\x$, is one for which, for all unit norm vectors $\z$, $\x' \z$ is sub-Gaussian. Also, its sub-Gaussian norm, $\|\x\|_{\varphi_2} = \sup_{\z:\|\z\|=1} \|\x' \z\|_{\varphi_2}$.

Let $K$ denote the maximum of the sub-Gaussian norms of the rows of $\A$. %Recall that $\A$ is a matrix with mutually independent, sub-Gaussian, and isotropic rows.
Theorem 5.39 of \cite{vershynin} shows that, for any $t>0$, with probability at least $1 - \exp(-c_Kt^2)$, the minimum singular value of $\A$ is more than $\sqrt{N} - (C_K\sqrt{n} + t)$ and the maximum is less than $\sqrt{N} + (C_K\sqrt{n} + t)$. Here $C_K$ and $c_K$ are numerical constants that depend only on $K$. These bounds are obtained by bounding the deviation of $\frac{1}{N} \A'\A$ from its expected value, which is equal to $\I$. %Remark 5.40 of \cite{vershynin} states a simple extension to the case when all rows of $\A$ have the same second moment matrix.

Our generalization of this result does two extra things. First, it bounds $\| \A'\A - \E[\A'\A]\|$, even when the different rows of $\A$ do not have the same second moment matrix.
Second, it states a separate result that bounds $\|\A\z\|^2$ for one specific vector $\z$. %$|\z'(\A'\A - \E[\A'\A]) \z |$ for a given vector $\z$.  are not isotropic and, in fact,
This bound clearly holds with much higher probability than the bound on $\|\A'\A - \E[\A'\A]\|$. The proof approach for getting our result is the same as that used to get Theorem 5.39 of \cite{vershynin}. Thus, our generalization would be obvious to a reader who understands the proof of that result. However, it is a useful addition to the literature for readers who would like to just use results from \cite{vershynin} in their work, without having to understand all their proof techniques.

%{\bf Notation. }
%The notation $\a_{i,k} \stackrel{\text{iid}}{\sim} \n(\mu, \Sigma)$ means that the vectors $\a_{i,k}$ are iid Gaussian with mean $\mu$ and covariance matrix $\Sigma$; and $\b_k \indepsim \n(\mu_k,\Sigma_k)$ means that the $\b_k$'s are mutually independent and $\b_k$ is generated from $\n(\mu_k,\Sigma_k)$.

%f?? We use $\|.\|_F$ to denote the Frobenius norm. (for all $i$ and $k$)
%The notation $\mathbbm{1}_{\zeta}$ denotes the indicator function for statement $\zeta$, i.e., $\mathbbm{1}_{\zeta}=1$ if $\zeta$ is true and $\mathbbm{1}_{\zeta}=0$ otherwise.
%For a vector $\z$, $|\z|$, $\sqrt{\z}$ and $\mathrm{phase}(\z)$ compute the {\em element-wise} magnitude, square-root, and phase\footnote{Here, and throughout, the term ``phase" actually refers to $e^{\sqrt{-1} \phi}$ where $\phi$ is the inverse tangent of the ratio between the imaginary and real parts of a complex number. Thus, for a real number the phase is $+1$ or $-1$.}
%We use `eigenvector' to always refer to unit norm eigenvector.
%Thus the output is a vector of the same size as $\z$.
%? remove We use $\e_i$ to denote the $i$-th column of the identity matrix $\I$.
%For a matrix $\bm{H}$, $\bm{H} \evdeq \bm{U\Lambda U}'$ denotes its reduced eigenvalue decomposition with the eigenvalues in $\Lambar$ arranged in non-increasing order. For Hermitian matrices

\section{Our Result}

\begin{theorem} \label{versh}
Suppose that $\w_j$, $j=1,2,\dots, N$, are $n$-length independent, sub-Gaussian random vectors with sub-Gaussian norms bounded by $K$. Let
\[
\bm{D}:=\frac{1}{N} \sum_{j=1}^N \w_j \w_j{}' - \frac{1}{N} \sum_{j=1}^N \E[ \w_j \w_j{}'].
\]
For an $\varepsilon > 0$,
\begin{enumerate}
%Thus, the probability of the above event is at least $ 1- 2 \exp(-cn)$ if $N > (c \log 9) n / \epsilon^2$.
\item for a given vector $\z$, with probability (w.p.) at least $1 - 2 \exp(-c \min(\varepsilon, \varepsilon^2) N)$,
\[
| \z' \bm{D} \z |  \le 4\varepsilon K^2 \| \z \|^2;
\]

\item w.p. at least $ 1 - 2 \exp(n \log 9 -c \min(\varepsilon, \varepsilon^2) N)$,
\[
\| \bm{D} \| \le  4\varepsilon K^2.
\]
\end{enumerate}
for a numerical constant $c$. %The constant $c$ is re-used multiple times for different numerical values.
\end{theorem}
{\em  Here, and throughout the paper, the letter $c$ is reused to denote different numerical constants.}

\begin{proof} The proof is given in the Supplementary Material. It follows using the approach developed in \cite{vershynin}.
%follows using the approach developed in \cite{vershynin} to prove its Theorem 5.39. We have included the proof of Theorem \ref{versh} in the Supplementary Material.
\end{proof}

\begin{remark} \em
Recall that $\A:= [\w_1, \w_2, \dots \w_N]'$ and so $\sum_j \w_j \w_j{}' = \A'\A$. Thus the first claim implies that $\frac{1}{N}\|\A \z \|^2 = \frac{1}{N} \z' \A' \A \z$ lies in the interval $(\frac{1}{N} \z' \E[ \A' \A] \z) \pm 4\varepsilon K^2  \| \z \|^2$ w.p. $\ge 1 - 2 \exp(-c \min(\varepsilon, \varepsilon^2) N)$.
Using Weyl's inequality, the second claim implies implies that w.p. $\ge 1 - 2 \exp(n \log 9 -c  \min(\varepsilon, \varepsilon^2) N)$, (a) $ \lambda_{\max}(\frac{1}{N} \A' \A)$ is smaller than $\lambda_{\max}(\frac{1}{N} \E[\A' \A]) + 4\varepsilon K^2$ and (b) $\lambda_{\min}(\frac{1}{N} \A' \A)$ is larger than $\lambda_{\min}(\frac{1}{N} \E[\A' \A]) - 4\varepsilon K^2$.
%Thus, with the same probability, $\|\A\|$ lies in $\sqrt{N}$ times this interval.
\end{remark}

%Theorem \ref{versh} is different from \cite[Theorem 5.39 or Remark 5.40]{vershynin} in 2 ways.%
%\begin{enumerate}
%\item It does not require $\w_j$'s to to be isotropic or to even have the same second moment matrix, $\E[\w_j \w_j{}']$.
%%be identically distributed. We use this to bound $\|\Y^- - \E[\Y^-|\B]\|$ conditioned on $\B$ in Lemma \ref{Y_minus}; here $\Y^-$ is a matrix that will be shown to be a close approximation of $\Y^U$. Conditioned on $\B$, observe that the rank one matrices summed to define $\Y^U$ are not identically distributed; the same is true for those used to define $\Y^-$.
%%This is important because, conditioned on $\B$, the r.v.'s
%%it is needed to bound $\|\Y^- - \E[\Y^-|\B]\|$ conditioned on $\B$ in Lemma \ref{Y_minus}; here $\Y^-$ is a matrix that will be shown to be a close approximation of $\Y^U$.
%\item  It states two separate results. Its first claim, which is true only for a specific vector $\z$, holds with a much higher probability than the second one.
%%Using the first claim in places where the second one is not needed is critical in obtaining a smaller lower bound on $m$ than $cn$.%what TWF needs.
%\end{enumerate}
%%Using the first claim where it can be applied is critical a lower bound on $m$  that is smaller than what TWF needs which is $m \ge c n$. $\|\Y_U - \E[\Y_U|\B]\|$ conditioned on $\B$ (we will actually bound

Theorem 5.39 of \cite{vershynin} is a corollary of the second claim of Theorem \ref{versh} specialized to isotropic $\w_j$'s. In that case $\E[\A' \A] = N \I$ and thus, by using the remark above with  $\varepsilon$ appropriately set (use $\varepsilon  = \frac{(\log 9) \sqrt{n}}{\sqrt{2} c\sqrt{N}} + \frac{t}{4K^2 \sqrt{N}}$), we get Theorem 5.39 of \cite{vershynin}.
%that $\|\A\|$ lies in the interval $\sqrt{N} \pm (4K^2 c \sqrt{n} + t)$ w.p. at least $1 - 2 \exp(n \log 9 -c \varepsilon^2 N) \ge 1 - 2 \exp(n \log 9 - n \log 9 - \frac{2t^2}{16 K^4}) ) = 1 - 2 \exp(-\frac{t^2}{8K^4})$. The last inequality used the fact that $(a+b)^2 \le 2a^2 + 2b^2$.

%It follows by letting the matrix $\A = [\w_1, \w_2, \dots \w_N]'$, setting $\varepsilon = \frac{c\sqrt{n} + t}{\sqrt{N}}$ and using Lemma 5.36 of \cite{vershynin}.

%Remark 5.40 more generally to $\w_j$'s with the same second moment matrix, i.e. for $\w_j$'s with $\E[\w_j \w_j{}'] = \bm\Sigma$

\subsection{An example application of Theorem \ref{versh}}
One example where both claims of the above result are useful but the result of Theorem 5.39 (or of Remark 5.40) of \cite{vershynin} does not suffice is in analyzing the initialization step of our recently proposed low-rank phase retrieval algorithm \cite{long_pr,icassp_pr}. In fact this is where we first used this generalization. %Without introducing the entire low rank phase retrieval problem, we use two examples motivated by the proof of the result of \cite{long_pr}
The example given below is motivated by this application. %e appl the proof of the results of \cite{long_pr}.
%We give next one example application of both claims of Theorem \ref{versh} next. Both these are either taken from or simplified versions of what was done in the proof of the results of \cite{long_pr}.

Consider $n$-length independent and identically distributed, standard Gaussian random vectors $\a_{i,k}$, i.e., $\a_{i,k} \iidsim \n(0,\I)$, with $i=1,2,\dots,m$ and $k=1,2,\dots,q$; and $n$-length deterministic vectors $\x_k$, $k=1,2,\dots,q$. Assume that $q \le n^2$.
Consider bounding
\[
b:=\left|\frac{1}{m} \sum_{i=1}^m (\a_{i,k}{}'\x_k)^2 - \|\x_k\|^2 \right|
\]
By applying item 1 of Theorem \ref{versh} with $N=m$ and $\w_j \equiv \a_{j,k}$, $b_k \le \epsilon \|\x_k\|^2$ w.p. at least $1- \exp(-c\epsilon^2 m)$. Such a bound holds for all $k=1,2,\dots,q$ w.p. at least $1- q \exp(-c\epsilon^2 m)$. Thus, to ensure that this bound holds w.p. at least $1 - 1/{\mathrm{poly}(n)}$, we need $m \ge \frac{ c (\log n + \log q) }{ \epsilon^2 }= \frac{ c \log n }{ \epsilon^2 }$ since $q \le n^2$. Here $\mathrm{poly}(n)$ means polynomial in $n$.

On the other hand,  to apply \cite[Theorem 5.39]{vershynin}, we first need to upper bound the $b_k$'s as
\begin{align*}
b_k & = |\x_k{}' \left(\frac{1}{m} \sum_{i=1}^m \a_{i,k} \a_{i,k}{}' - \I \right) \x_k|  \\
 & \le \|\x_k\|^2 \left\|\frac{1}{m} \sum_{i=1}^m \a_{i,k} \a_{i,k}{}' - \I \right\|
\end{align*}
With this, we can get the same bound as above on the $b_k$'s by applying \cite[Theorem 5.39]{vershynin} with $ t = \sqrt{m} 4K^2 \epsilon - C_K \sqrt{n}$ (or, equivalently, by applying item 2 of Theorem \ref{versh} above). But the bound would hold with probability lower bounded by $1- \exp(n \log 9 - c\epsilon^2 m)$. For a given $m$, this is a much smaller probability. Said another way, one would need $m \ge \frac{c n}{ \epsilon^2}$ for the probability to be high enough (at least $1 - 1/{\mathrm{poly}(n)}$). This is a much larger lower bound on $m$ than the earlier one.
%to be at least $1 - \frac{1}{\mathrm{poly}(n)}$. Of course with $m \ge cn$, the result would hold w.p. $1-  \exp(-cn)$ which is actually even higher.

To see an application of item 2 of Theorem \ref{versh}, consider bounding
\[
\tilde{b}:= \left\|\frac{1}{mq} \sum_{k=1}^q \sum_{i=1}^m \a_{i,k} \a_{i,k}{}' f_k^2 - \frac{1}{q}\sum_{k=1}^q f_k^2 \right\|
\]
where $f_k$'s are scalars. By conditioning on the $f_k$'s, we can apply item 2 of Theorem \ref{versh} on all the $N=mq$ vectors $(\a_{i,k} f_k)$ to conclude that $\tilde{b} \le \epsilon_2 \max_k f_k^2$, w.p. at least $ 1- 2 \exp(n \log 9 - c \epsilon_2^2 mq)$. %Notice that we are able to apply the result with $N = mq$.
Thus, the bound holds w.p. at least $1 - 1/{\mathrm{poly}(n)}$ if $m \ge  \frac{c n}{q \epsilon_2^2}$.

Observe that the $\a_{i,k}$'s are isotropic independent sub-Gaussian vectors but $\a_{i,k} f_k$'s are not. In fact, $\E[\a_{i,k} \a_{i,k}{}' f_k^2] = f_k^2$ and hence the vectors $\a_{i,k} f_k$ also do not have the same second moment matrix for all $k,i$. As a result, we cannot apply Theorem 5.39 or Remark 5.40 of \cite{vershynin} to bound $\tilde{b}$ if we want to average over all the $mq$ vectors. To apply one of these, we first need to upper bound $\tilde{b}$ as
\[
\tilde{b} \le \frac{1}{q} \sum_{k=1}^q \left\|\frac{1}{m} \sum_{i=1}^m \a_{i,k} \a_{i,k}{}' - \I \right\| f_k^2
\]
Now using \cite[Theorem 5.39]{vershynin}, we get $\tilde{b} \le \epsilon_2 \frac{1}{q}  \sum_{k=1}^q  f_k^2 \le \epsilon_2  \max_k f_k^2 $ w.p. at least $ 1- 2 \exp(n \log 9 - c \epsilon_2^2 m)$. Observe that $mq$ is replaced by $m$ in the probability now. Thus, to get the probability to be high enough (at least $1 - \frac{1}{\mathrm{poly}(n)}$) we will need $m \ge \frac{c n}{ \epsilon_2^2}$ which is, once again, a much larger lower bound than what we got by applying item 2 of Theorem \ref{versh}.

To understand the context, in \cite{long_pr}, $m$ is the sample complexity required for the initialization step of low-rank phase retrieval to get an estimate of the low-rank matrix $\X:=[\x_1, \x_2, \dots \x_q]$ that is within a relative error $c\epsilon$ of the true $\X$ with probability at least $1-1/\mathrm{poly}(n)$. If we directly use the result from \cite{vershynin}, we will need $m \ge cn/\epsilon^2$, where as if we use Theorem \ref{versh}, we can get a lower bound that is smaller than $cn$ (when $q$ is large enough).

%Observe that $\a_{i,k}$'s are sub-Gaussian with sub-Gaussian norm bounded by a constant.

%To see why the above result is more useful than Theorem 5.39 or Rema of \cite{vershynin}

\section{Conclusions}
We proved a simple generalization of a result of Vershynin \cite{vershynin} for random matrices with independent, sub-Gaussian rows. %In a similar fashion, it should be possible to also obtain   generalizations of some of the other results proved in \cite{vershynin} too.

We should mention that the first claim of Theorem \ref{versh} can be further generalized for two different vectors $\z_1$ and $\z_2$ as follows: with the same probability, $|\z_1{}' \bm{D} \z_2| \le 4 \varepsilon^2 K^2 (\|\z_1\|^2 + \|\z_2\|^2)$. This follows because, for two sub-Gaussian scalars, $x, y$, $xy$ is sub-exponential with sub-exponential norm bounded by $c(\|x\|_{\psi_2}^2+\|y\|_{\psi_2}^2)$ \cite{vince_vu_nips}.
%, for example, Theorem 5.41 which holds for matrices with independent, heavy-tailed rows.
%, Theorem 5.39 of

\bibliographystyle{IEEEbib}
%\bibliography{../../bib/tipnewpfmt_kfcsfullpap}
\bibliography{tipnewpfmt_kfcsfullpap}
\appendix
\section*{Proof of our result}

%Based on reviewers' comments, this Appendix could either be supplementary material or it could be part of the paper.

\subsection{Preliminaries} \label{prelim}

%\subsubsection{$\epsilon$-net}
%All material in this and the next section is taken from \cite{vershynin}. \cite[Definition 5.1]{vershynin}}
As explained in \cite{vershynin}, nets are a convenient means to discretize compact metric spaces. The following definition is \cite[Definition 5.1]{vershynin} for the unit sphere.
%Here we will need to discretize the unit sphere in $\Re^n$. Hence we define a net only for that specific case. The distance metric will be the Euclidean norm distance. %, denoted $\S^{n-1}$
%\begin{definition}[$\epsilon$-net and covering number of the unit sphere in $\Re^n$]
%We define a net only on the unit sphere in $\Re^n$, denoted $\S^{n-1}$. The metric is  the usual Eucliden norm metric.
%
For an $\epsilon > 0$, a subset $\N_\epsilon$ of the unit sphere in $\Re^n$ is called an $\epsilon$-net if, for every vector $\x$ in the unit sphere, there exists a vector $\y \in \N_\epsilon$ such that $\|\y - \x\| \le \epsilon$.

The covering number of  the unit sphere in $\Re^n$, is the minimal cardinality of an $\epsilon$-net on it. In other words, it is the size of the smallest $\epsilon$-net, $\N_\epsilon$, on it.
%\end{definition}%the unit sphere  in $\Re^n$,

\begin{fact} \em \label{fact_nets}
\ %minimum cardinality of an $\epsilon$-net on    (its covering number)
\begin{enumerate}
\item By Lemma 5.2 of \cite{vershynin}, the covering number of the unit sphere in $\Re^n$ is upper bounded by $(1+\frac{2}{\epsilon})^n$. %Thus, for example, if $\epsilon=1/4$, then the covering number is at most $9^n$.

\item By Lemma 5.4 of \cite{vershynin}, for a symmetric matrix, $\W$,
$
\|\W\| \le \max_{\x: \|\x\|=1} \|\x' \W \x\| \le \frac{1}{1-2 \epsilon} \max_{\x \in \N_{\epsilon} } \|\x' \W \x\|.
$
\end{enumerate}
Thus, if $\epsilon = 1/4$, then $\|\W\| \le 2 \max_{\x \in \N_{1/4}  } \|\x' \W \x \|$ and the cardinality of the smallest such net is at most $9^n$.
\end{fact}

%\subsubsection{sub-Gaussian and sub-exponential r.v.s}
%A r.v. $x$ is sub-Gaussian if and only if the following holds: there exists a constant $K$ such that
%$
%\E[|x|^p]^{1/p} \le K \sqrt{p}
%$
%for all integers $p \ge 1$.
%The smallest such $K$ is referred to as the sub-Gaussian norm of $x$, denoted $\|x\|_{\varphi_2}$. Thus,
%$
%\|x\|_{\varphi_2} = \sup_{p \ge 1} p^{-1/2} \E[|x|^p]^{1/p}.
%$
%
%A random vector $\x$ is sub-Gaussian if, for all unit norm vectors $\z$, $\x' \z$ is sub-Gaussian. Also, $\|\x\|_{\varphi_2} = \sup_{\z:\|\z\|=1} \|\x' \z\|_{\varphi_2}$.

A r.v. $x$ is sub-exponential if the following holds: there exists a constant $K_e$ such that
$
\E[|x|^p]^{1/p} \le K_e {p}
$
for all integers $p \ge 1$; the smallest such $K_e$ is referred to as the sub-exponential norm of $x$, denoted  $\|x\|_{\varphi_1}$ \cite[Section 5.2]{vershynin}. %Thus, There are other equivalent definitions, see \cite[Section 5.2]{vershynin}.
%\[
%\|x\|_{\varphi_1} = \sup_{p \ge 1} p \E[|x|^p]^{1/p}.
%\]
%There are other equivalent definitions, see \cite[Section 5.2]{vershynin}.There are other equivalent definitions, see \cite[Section 5.2]{vershynin}.

The following facts will be used in our proof.% in our proof. They follow using lemmas in \cite{vershynin}.
\begin{fact} \em \label{fact_versh} \
\begin{enumerate}
\item  \label{sub_expo}
If $\x$ is a sub-Gaussian random vector with sub-Gaussian norm $K$, then for any vector $\z$, (i) $\x'\z$ is sub-Gaussian with sub-Gaussian norm bounded by $K \|\z\|$; (ii) $(\x'\z)^2$ is sub-exponential with sub-exponential norm bounded by $2 K^2 \|\z\|^2$; and (iii) $(\x'\z)^2 - \E[(\x'\z)^2]$ is centered (zero-mean), sub-exponential with sub-exponential norm bounded by  $4 K^2 \|\z\|^2$. This follows from the definition of a sub-Gaussian random vector; Lemma 5.14 and Remark 5.18 of \cite{vershynin}.

\item \label{cor_5_17}
By \cite[Corollary 5.17]{vershynin}, if $x_i$, $i=1,2, \dots N$, are a set of independent, centered, sub-exponential r.v.'s with sub-exponential norm bounded by $K_e$, then,
%for an $\epsilon>0$ that satisfies $\frac{\epsilon}{K_i} < 1$,  $\Pr( | \frac{1}{N} \sum_i Z_i| > \epsilon) \le 2 \exp\left(c \frac{\epsilon^2}{K^2} N \right)$.  Equivalently,
for any $\varepsilon >0$,
\[
\Pr\left( |  \sum_{i=1}^N x_i| > \varepsilon K_e N \right) \le 2 \exp( - c \min(\varepsilon, \varepsilon^2)  N).
\]

\item \label{K_gauss_vec}
If $\x \sim \n(0, \Lambar)$ with $\Lambar$ diagonal, then $\x$ is sub-Gaussian with $\|\x\|_{\varphi_2} \le c \sqrt\lammax$.
%Moreover, if  $\y = [x_1, \x']'$ where $x_1$ is a  zero mean bounded r.v. with bound $M$ and is independent of the vector $\x$, then $\|\y\|_{\varphi_2} \le c \max(M, \lammax)$. This follows from Example 5.8 and Lemma 5.24 of \cite{vershynin}.

%\item \label{rayleigh_type_bnd} If $\x_i \sim \n(0, \Lambar)$, for $i=1,2,\dots, N$, are $n$-length random vectors and  $\Lambar$ is diagonal, then
%\begin{align*}
%& \Pr \left( \max_{i=1,2,\dots, N} \|\x_i\|^2 \le \lammax \cdot n \cdot 2 \nu \right)   \\
%& \ge 1 - 2n N \exp(-\nu) , \ \text{for} \   \nu > 1.
%\end{align*}
%This is a direct consequence of eq. 5.5 of \cite{vershynin} which says that if $x \sim \n(0,1)$, then $\Pr(|x_i| > t) \le 2 \exp(-t^2/2)$ for a $t>1$.
%Using this along with the union bound first for bounding $\|\x_i\|^2 = \sum_{j=1}^n (\x_i)_j^2$ for a given $i$ and then for bounding its $\max$ over $i$ gives the above result.
%
%\item \label{subG_bnd}
%Using \cite[Lemma 5.5]{vershynin}, if $\x_i$'s are sub-Gaussian random vectors with sub-Gaussian norm bounded by $K$, then the following generalization of the above fact holds:
%$\Pr \left( \max_{i=1,2,\dots, N} \|\x_i\|^2 \le K^2 \cdot n \cdot 2 \nu \right)  \ge 1 - Cn N \exp(- c\nu)$.% This follows using \cite[Lemma 5.5]{vershynin} and the union bound.
%%, \ \text{for} \   \nu > 1.

\end{enumerate}
\end{fact}

\subsection{Proof of Theorem \ref{versh}} \label{proof_versh}

%\begin{proof}[Proof of Theorem \ref{versh}]
%The first part follows using Corollary 5.17 of \cite{vershynin} (for sums of sub-exponential random variables). The second part uses  an $\epsilon$-net argument followed by application of Corollary 5.17 of \cite{vershynin}. $The first part holds with a h
% (for sums of sub-exponential random variables). The proof of the first part follows without use of the $\epsilon$-net argument. %The complete proof is provided in Appendix \ref{proof_versh}.
%
%We prove the first part first.
The proof strategy is similar to that of Theorem 5.39 of \cite{vershynin}.
By Fact \ref{fact_versh}, item \ref{sub_expo}, for each $j$, the r.v.s $\w_j{}' \z$ are sub-Gaussian with sub-Gaussian norm bounded by $K \| \z\|$;  $(\w_j{}' \z)^2$ are sub-exponential with sub-exponential norm bounded by $2K^2 \| \z\|^2$; and $(\w_j{}' \z)^2 - \E[(\w_j{}' \z)^2] = \z'(\w_j \w_j') \z - \z'(\E[\w_j \w_j']) \z$ are centered sub-exponential with sub-exponential norm bounded by $4K^2 \| \z\|^2$. Also, for different $j$'s, these are clearly mutually independent. Thus, by applying Fact \ref{fact_versh}, item \ref{cor_5_17} (Corollary 5.17 of \cite{vershynin}) with $K_e = 4K^2 \| \z\|^2$ we get the first part.

To prove the second part, let $\N_{1/4}$ denote a 1/4-th net on the unit sphere in $\Re^n$. Let $\W:= \frac{1}{N} \sum_{j=1}^N ( \w_j \w_j' - \E[ \w_j \w_j'])$.
Then by Fact \ref{fact_nets} (Lemma 5.4 of \cite{vershynin})
\beq \label{epsnet_bnd}
\|\W\| \le 2 \max_{\z \in  \N_{1/4}} |\z' \W \z|
\eeq
Since $\N_{1/4}$ is a finite set of vectors, all we need to do now is to bound $|\z' \W \z|$ for a given vector $\z$ followed by applying the union bound to bound its  maximum over all ${\z \in  \N_{1/4}}$.
The former has already been done in the first part.
%For a given vector $\z$, $|\z' \W \z|$ has already been bounded in the first part.
By Fact \ref{fact_nets} (Lemma 5.2 of \cite{vershynin}), the cardinality of $\N_{1/4}$ is at most $9^n$. Thus, using the first part, %and taking the union bound over all vectors $\z \in \N_{1/4}$,
$
\Pr\left(\max_{\z \in  \N_{1/4}} |\z' \W \z|  \ge  \frac{4\varepsilon K^2}{2} \right) \le 9^n \cdot 2 \exp(-c \frac{\min(\varepsilon, \varepsilon^2)}{4} N) =  2 \exp(n \log 9 -c \varepsilon^2 N).
$
By \eqref{epsnet_bnd}, we get the result.
%\end{proof}

\end{document}